\documentclass[12pt]{amsart}
\usepackage{amscd,amssymb}
\usepackage[arrow,matrix]{xy}
\usepackage[plainpages,backref,urlcolor=blue]{hyperref}

\theoremstyle{plain}
\newtheorem{thm}[subsection]{Theorem}
\newtheorem{lem}[subsection]{Lemma}
\newtheorem{prop}[subsection]{Proposition}
\newtheorem{cor}[subsection]{Corollary}

\theoremstyle{definition}
\newtheorem{rk}[subsection]{Remark}
\newtheorem{definition}[subsection]{Definition}
\newtheorem{ex}[subsection]{Example}

\numberwithin{equation}{section}
\newcommand{\OO}{{\mathcal O}}

\newcommand{\DD}{{\mathcal D}}
\newcommand{\CC}{{\mathcal C}}
\newcommand{\LL}{{\mathcal L}}
\newcommand{\I}{{\mathcal I}}

\newcommand{\be}{{\beta}}
\newcommand{\NN}{{\mathcal N}}

\newcommand{\N}{\mathbb{N}}
\newcommand{\Z}{\mathbb{Z}}

\newcommand{\C}{\mathbb{C}}

\newcommand{\PP}{\mathbb{P}}
\newcommand{\HH}{\mathbb{H}}

\DeclareMathOperator{\defect}{def}
\DeclareMathOperator{\codim}{codim}

\begin{document}

\title [On the syzygies and Hodge theory of  nodal hypersurfaces]
{On the syzygies and Hodge theory of nodal hypersurfaces }

\author[Alexandru Dimca]{Alexandru Dimca$^1$}
\address{Universit\'e C\^ ote d'Azur, CNRS, LJAD, France }
\email{dimca@unice.fr}

\thanks{$^1$ Partially supported by Institut Universitaire de France.}

\subjclass[2000]{Primary  14C30, 13D40; Secondary  32S35, 13D02}

\keywords{nodal hypersurfaces, syzygies, mixed Hodge structure, pole order filtration}

\begin{abstract} We give sharp lower bounds for the degree of the syzygies involving the partial derivatives of a homogeneous polynomial defining an even dimensional nodal hypersurface. 
This implies the validity of formulas due to M. Saito, L. Wotzlaw and the author for the graded pieces with respect to the Hodge filtration of the top cohomology of the hypersurface complement in many new cases. A classical result by Severi on the position of the singularities of a nodal surface in $\PP^3$ is improved and  applications to deformation theory of nodal surfaces are given.

\end{abstract}

\maketitle


\noindent {\it Dedicated to the memory of Alexandru Lascu, who was always searching for the Truth, and encouraging  others to do the same.}

\section{Introduction} \label{sec:intro}

Let $S=\C[x_0,...,x_n]$ be the graded ring of polynomials in $x_0,,...,x_n$ with complex coefficients and denote by $S_r$ the vector space of homogeneous polynomials in $S$ of degree $r$. 
For any polynomial $f \in S_r$, we define the {\it Jacobian ideal} $J_f \subset S$ as the ideal spanned by the partial derivatives $f_0,...,f_n$ of $f$ with respect to $x_0,...,x_n$ respectively, 
and we define the corresponding graded {\it Milnor} (or {\it Jacobian}) {\it algebra} by
\begin{equation} 
\label{eq1}
M(f)=S/J_f.
\end{equation}

The Milnor algebra $M(f)$ can be seen (up to a twist in grading) as the top cohomology $H^{n+1}(K^*(f))$,
 where $K^*(f)$ is the Koszul complex of $f_0,...,f_n$ with the natural grading $|x_j|=|dx_j|=1$, defined by 
\begin{equation} 
\label{Koszul}
K^*(f) : \  \   \    \   \   \     0 \to \Omega^0 \to \Omega^1 \to ...   \to \Omega^{n+1}\to 0
\end{equation}
with all the arrows  given by the wedge product by $df=f_0dx_0+f_1dx_1+...+f_ndx_n$.
The homogeneous components of the next cohomology group, say $H^n(K^*(f))_{n+r}$, describe the syzygies
$$\sum_{j=0,n}a_jf_j=0$$
where $a_j \in S_r$, modulo the trivial, or Koszul, syzygies generated by 
$$(f_j)f_i+(-f_i)f_j=0$$
for all $i<j$.
The following result was proved in \cite{DSt3}.

\begin{thm}
\label{thmA} Let $D:f=0$ be a nodal hypersurface of degree $d>2$ in $\PP^n$.

\smallskip

\noindent (i) If $n=2n_1+1$ is odd, then $H^n(K^*(f))_{m}=0$ for any $m \leq n_1d$.

\smallskip

\noindent (ii) If $n=2n_1$ is even, then $H^n(K^*(f))_{m}=0$ for any $m \leq n_1d-1$.

\end{thm}

Examples involving Chebyshev hypersurfaces shows that the bounds are best possible for $n$ even,
but not for $n$ odd, see \cite{DSt3}. The first purpose of this paper is to establish the following new bound in the case  $n$ odd, which is optimal, see Example \ref{ex1}.

\begin{thm}
\label{thmB} Let $D:f=0$ be a nodal hypersurface of degree $d>2$ in $\PP^n$.

If $n=2n_1+1$ is odd, then $H^n(K^*(f))_{m}=0$ for any $m \leq (n_1+1)d-\left[ \frac{d}{2} \right]-1$.

\end{thm}

A classical result of Severi \cite{Sev} says that if a surface $D$ in $\PP^3$ of degree $d$ has only nodes as singularities, then the set of nodes imposes independent conditions on hypersurfaces of degree $2d-5$. A modern proof of this result can be found in \cite{Nobile}, and a more general version is stated in \cite[Corollary H]{MP}. Using Theorem \ref{thmB} we can prove the following stronger version of Severi's result.

\begin{cor} 
\label{corSeveri}
Let $D:f=0$ be a nodal surface of degree $d>2$ in $\PP^3$. Then the set of nodes imposes 
independent conditions on hypersurfaces of degree 
$$d+\left[ \frac{d}{2} \right]-3.$$
\end{cor}

The surfaces in $\PP^3$ of degree $d$ are parametrized by $\PP^N$, with $N={d+3 \choose 3}-1$. There is a locally closed subscheme $V_{d,n} \subset \PP^N$ parametrizing the surfaces $D \subset \PP^3$ of degree $d$ and having exactly $n$ nodes, see \cite{Nobile}, \cite{Se0}. Let $X$ be the minimal resolution of a surface $D$  corresponding to a point in 
$V_{d,n}$ and let $R(X)$ be the formal moduli space of $X$ as in \cite{BW}, regarded as a 
complete local ring.
The first two claims of the following result are given in \cite[Theorem 3.2]{Nobile} for the case $d=5$. The similar results in the case of nodal plane curves are classical, see 
\cite{Nobile0} and \cite[Corollary 4.7.8, Theorem 4.7.18 and Corollary 4.7.19]{Se0}.

\begin{thm} 
\label{corNobile} Fix a degree $d \in \{5,6,7\}$ and a positive integer $n$.

\begin{enumerate} 

\item  The variety $V_{d,n}$ is a smooth locally closed subscheme of $\PP^N$. If $V_{d,n}\ne \emptyset$, then $V_{d,n}$ has pure dimension $N-n$.

\item If $V_{d,n}\ne \emptyset$, then $V_{d,n'}\ne \emptyset$ for any positive integer $n' \leq n$.

\item The surface $X$ is unobstructed, i.e. the complete local ring $R(X)$ is regular, for any nodal surface $D \in V_{d,n}$.

\end{enumerate} 

\end{thm}

In their paper \cite{BW}, the authors give several examples of obstructed nodal surfaces for any degree $d \geq 8$. Our Theorem \ref{corNobile} explains why there are no such examples for $d<8$.

The idea of proof of Theorem \ref{thmB} is similar to that used in the proof of 
Theorem \ref{thmA}, namely the interplay between Hodge filtrations, pole order filtrations
and some spectral sequences. In the case of Theorem \ref{thmA} it was enough to look at the cohomology of the complement $U=\PP^n \setminus D$, while in the case at hand we have to consider the eigenspaces $H^*(F)_{\lambda}$ of the monodromy action on the corresponding Milnor fiber $F: f(x)-1=0$, which is technically more complicated.

The second aim of this paper is to show that the information we have obtained on the syzygies
allows us to prove that an algebraic description of some of the graded pieces $Gr_F^pH^n(U,\C)$ of the top cohomology group $H^n(U,\C)$ with respect to the Hodge filtration $F^*$
given in \cite{DSW} holds in many cases, see Theorem \ref{thmC}. For basic facts on mixed Hodge structures we refer to the excellent monograph \cite{PS}.

This paper was written back in 2013, and some of the results here were improved, both in generality and in presentation, in our subsequent joint work with Morihiko Saito, see \cite{DS4}, \cite{DS5}. For instance, Theorem \ref{thmB} is generalized to hypersurfaces with isolated weighted homogeneous singularities in \cite[Theorem 9]{DS5}. However,  these preprints use a different approach and are not in the final form. Since our results have been already quoted in important papers such as \cite{MP}, we have decided to publish them
as an alternative view-point on this subject. 

On the other hand, Corollary \ref{corSeveri} and Theorem \ref{corNobile} above are new.

\medskip

The author would like to thank Morihiko Saito and Edoardo Sernesi for very useful discussions related to this paper.

\section{Spectral sequences for Milnor fibers of homogeneous polynomials} \label{sec2}

For  a fixed integer $k=1,...,d$, we set $\lambda =\exp(-2 \pi ik/d)$ and let $L_k$ be the rank one local system on $U$ such that
\begin{equation} 
\label{eq2}
H^*(U,L_k)=H^*(F)_{\lambda},
\end{equation}
see for details \cite{D2}, p. 211 and note that here the eigenspaces are with respect to the local system monodromy $T$  acting on $H^*(F)$, as explained in \cite{DS3} . 

Let $j:U \to X=\PP^n$ be the inclusion, and let $\LL_k$ be the meromorphic extension of $L_k \otimes_{\C}\OO_U$. Then $\LL_k$ is a regular holonomic $\DD_X$-module,
see \cite{Sa2}, section (4.8), and one clearly has
\begin{equation} 
\label{eq3}
\HH^*(X,DR(\LL_k))=\HH^*(X,Rj_*L_k)=H^*(F)_{\lambda}.
\end{equation}
The $\DD_X$-module $\LL_k$ has a natural (increasing) pole order filtration $P_*$ such that 
$P_j\LL_k$ is isomorphic to $\OO_X((jd+k)$ for $j \in \N$ and $P_j\LL_k=0$ for $j<0$, see
 (3.1.3) in \cite{DS4}. We define a decreasing filtration $P^*$ on $\LL_k$ by putting $P^m=P_{-m}$ for any $m \in \Z$.
 Then the de Rham complex $DR(\LL_k)$ has an induced decreasing filtration
$P^*$ and this induces the {\it pole order filtration} $P^*$ on the eigenspaces $H^*(F)_{\lambda}$
of the Milnor fiber cohomology. One has the following fundamental inclusion
\begin{equation} 
\label{HP1}
F^sH^j(F)_{\lambda} \subset P^sH^j(F)_{\lambda},
\end{equation}
for any $s$ and any $j$, where $F^s$ denotes the canonical Hodge filtration on the cohomology of the smooth quasi-projective variety $F$ constructed by P. Deligne, see (3.1.3), (4.4.7) and (4.4.8) in \cite{DS4}.  This extends the similar result for $H^*(U)$, see \cite{DD} and \cite{Sa1}, which correspond to the case $\lambda=1$.

On the other hand, for each $k=1,...,d$ there is a spectral sequence
\begin{equation} 
\label{sp2}
E^{p,q}_{1}(f,k)=H^{p+q}(K^*(f))_{qd+k} \Rightarrow H^{p+q-1}(F)_{\lambda}
\end{equation}
coming from the graded Gauss-Manin complex $C^*_f$ associated with $f$, see (4.4.4) and (4.5.3) in  \cite{DS4}. A similar spectral sequence is obtained by using the algebraic microlocal Gauss-Manin complex $\tilde C^*_f$ and the corresponding limit is the reduced cohomology of the Milnor fiber, see the   (4.5) in \cite{DS4}. In fact, one has
\begin{equation} 
\label{sp3}
H^{j+1}( C^*_f)_{k}=H^{j+1}( \tilde C^*_f)_{k}= H^{j}(F)_{\lambda}
\end{equation}
for any $k=1,...,d$ and $j>0$, see (4.2.3)  in \cite{DS4}.
These spectral sequences induce a filtration $P'$ on $H^*(F)_{\lambda}$
and one has $P^s=P'^{s+1}$, see (4.4.7) in  \cite{DS4}. When the hypersurface $D:f=0$ has only weighted homogeneous singularities, these spectral sequences degenerate at $E_2$, see the very recent paper by Morihiko Saito \cite{Sa3}.

\begin{rk}
\label{rkMB}
Remark 4.11 in \cite{Sa2} gives a very explicit description of the pole order filtration on the top cohomology group $H^n(F)$. Assume that we have a (finite) family of monomials $(g_j(x))_{j \in J}$ in $S$ such that the cohomology classes $[\omega_j]$ of the differential forms $\omega_j=g_j(x)\Delta(dx_0 \wedge ...\wedge dx_n)$ for $j \in J$ yield a basis of the $\C$-vector space $H^n(F)$. Then $Gr_P^{n-q}H^n(F)_{\lambda}$ is spanned by the classes $[\omega_j]$    with $\deg \omega_j = qd+k.$ Note that the eigenspace $Gr_P^{n-q}H^n(F)_1$ is spanned exactly by the classes $[\omega_j]$ in $Gr_P^{n-q}H^n(F)$ having a maximal degree $\deg \omega_j =(q+1)d$. In particular, the pole order filtrations on $H^n(F)_1$ and on $H^n(U)$, the latter constructed using a generalization of Griffiths approach in  \cite{D1}, Chapter 6 or in  \cite{DSt2}
are the same, as it should be since we have a natural isomorphism $H^n(F)_1=H^n(U)$.
Note that the corresponding basis for $Gr_P^{n-q}H^n(U)$ is usually written as
$$\sigma_j=\omega_j/f^{q+1}.$$
Similarly, the limit term $E^{n-q+1,q}_{\infty}(f,k)$, which is isomorphic to $Gr_P^{n-q}H^n(F)_{\lambda}$, has a basis given with the above notation by
$$\eta_j=g_j(x)\cdot dx_0 \wedge ...\wedge dx_n.$$

\end{rk}

\section{Hodge filtration versus pole order filtration on Milnor fibers} \label{sec3}

If $D$ is a nodal hypersurface in $\PP^n$, then one has an equality
\begin{equation} 
\label{HP2}
F^sH^n(U)= P^sH^n(U),
\end{equation}
for any $s \geq n-m+1$, with $m=\alpha_D=\frac{n}{2}$, see Corollary (0.12) in M. Saito \cite{Sa1}, or  the formula (1.1.3) in \cite{DSW}. The purpose of this section is to prove the following
similar result for the associated Milnor fibers.

\begin{prop}
\label{propA}
Let $D:f=0$ be a nodal hypersurface of degree $d>2$ in $\PP^n$ with $n$ odd. Then 
$$F^sH^n(F)_{\lambda}= P^sH^n(F)_{\lambda},$$
for any $s \geq n-m+1$ and any $\lambda \in \mu_d$, $\lambda \ne 1$, with $m=\frac{n}{2}$.
\end{prop}

\proof
Consider the hypersurface $\tilde D$ defined in $\PP^{n+1}$ by the equation $\tilde f(x,t)=0$,
with $\tilde f(x,t)=f(x)-t^d$. Let $\tilde U=\PP^{n+1} \setminus \tilde D$ and let $H:t=0$ be the hyperplane at infinity in $\PP^{n+1}$ such that $\tilde U \cap H=U$. Consider the Gysin long exact sequence
\begin{equation} 
\label{G1}
...\to H^{n-1}(U) \to H^{n+1}(\tilde U) \to H^{n+1}(\tilde U \setminus U) \to H^n(U) \to...
\end{equation}
The group $\mu_d$ of $d$-roots of unity acts on $\PP^{n+1}$ via
$$ \be \cdot (x_0:...:x_n:t)=(x_0:...:x_n:\be t).$$
This action extends the  action of $\mu_d$ on $\C^{n+1}=\PP^{n+1} \setminus H$, which is used to define the local system monodromy  $T:F \to F$, namely
$$ \be \cdot (x_0,...,x_n)=(\be^{-1} x_0,...,\be^{-1}  x_n).$$

It follows that the Gysin exact sequence \eqref{G1} inherits a $\mu_d$-action, such that for any
$\lambda \in \mu_d$, $\lambda \ne 1$, one has the following isomorphism of eigenspaces
\begin{equation} 
\label{G2}
i^*: H^{n+1}(\tilde U)_{\lambda} \to H^{n+1}(\tilde U \setminus U)_{\lambda}.
\end{equation}
On the other hand, one has $\tilde U \setminus U=\C^{n+1} \setminus F$, and hence a new Gysin sequence
shows that one has an isomorphism
\begin{equation} 
\label{G3}
R:H^{n+1}(\tilde U \setminus U)_{\lambda} \to  H^n(F)_{\lambda},
\end{equation}
induced by a residue morphism $R$ which has Hodge type $(-1,-1)$. By composing the above two isomorphisms, we get isomorphisms
\begin{equation} 
\label{G4}
Ri^*: F^{s+1}H^{n+1}(\tilde U)_{\lambda} \to F^sH^n(F)_{\lambda},
\end{equation}
for any integer $s$ and any $\lambda \in \mu_d$, $\lambda \ne 1$.

Now we look at the corresponding $P^*$ filtrations and  show  that
\begin{equation} 
\label{G5}
\dim P^{s+1}H^{n+1}(\tilde U)_{\lambda} =\dim P^sH^n(F)_{\lambda},
\end{equation}
for any integer $s$ and any $\lambda \in \mu_d$, $\lambda \ne 1$.

The cohomology of the filtered algebraic microlocal Gauss-Manin complexes $\tilde C^*_f$ and $\tilde C^*_{\tilde f}$ are closely related, namely
$$H^{\ell+1}( \tilde C^*_{\tilde f},P'  )=H^{\ell}(  \tilde C^*_f ,P'  )\otimes H^1(\tilde C^*_{t^d},P')$$
see (4.9) in \cite{DS4}. Looking at the homogeneous components corresponding to $k=d$ and taking the $\lambda$-eigenspaces yields the equality \eqref{G5} in view of \eqref{sp3}.

Finally, any singularity of the hypersurface $\tilde D$ has type $A_{d-1}$, and hence the corresponding $\alpha_{\tilde D}$ is exactly $\tilde m=\frac{n}{2}+\frac{1}{d}$.
Using Corollary (0.12) in M. Saito \cite{Sa1}, or the formula (1.1.3) in \cite{DSW}, we see that
\begin{equation} 
\label{HP3}
F^sH^{n+1}(\tilde U)= P^sH^{n+1}(\tilde U),
\end{equation}
for any $s \geq n-\tilde m+2$. Using the formulas \eqref{G4} and \eqref{G5} and the inclusion \eqref{HP1} we complete the proof, as one may clearly replace $\tilde m$ by $m$ as soon as $d>2$ and $n$ is odd.
 
\endproof

\begin{rk}
\label{rkHvsP} Let $D:f=0$ be a   curve in $\PP^2$, with arbitrary isolated singularities.
In this case an algorithm computing the dimension of the eigenspaces $H^{m}(F,\C)_{\lambda }$, for $m=1,2$ and the dimension of the graded pieces of the pole order filtration on the Milnor fiber cohomology is described in \cite{DStgenMF}. Moreover, one has the following.

 (i) The  Hodge filtration $F$ and the pole order filtration $P$ coincide on
 $H^{1}(F,\C)_{\lambda }$ for any $\lambda$ in all the examples we have computed so far. We {\it conjecture} that the two filtrations $F$ and $P$ coincide on $H^{1}(F,\C)_{\lambda }$ always.
 
 (ii) On the other hand, the two filtrations $F$ and $P$ do not coincide on $H^{2}(F,\C)_{\lambda }$
even in very simple cases, e.g. $C:f=(x^2-y^2)(x^2-z^2)(y^2-z^2)=0$ and ${\lambda }=-1$.
A computation of the Hodge filtration on $H^{2}(F,\C)$ in this case can be found in \cite{PBthesis}. Note  also that the mixed Hodge structure on 
$H^2(F,\C)_{\ne 1}$ is not pure in general. For a line arrangement, one can use the formulas for the spectrum given in \cite{BS} to study the interplay between monodromy and Hodge filtration on 
$H^2(F,\C)_{\ne 1}$.

\end{rk}

\section{Proof of  Theorem \ref{thmB} and some examples} \label{sec4}

Assume from now on that $n=2n_1+1$. Then the equality of filtrations in Proposition \ref{propA}
holds for $s \geq n_1+2$. We show now that this equality may be extended as follows.

\begin{lem}
\label{L}
Let $D:f=0$ be a nodal hypersurface of degree $d>2$ in $\PP^{2n_1+1}$ and set $s=n_1+1$. Then 
$$F^sH^n(F)_{\lambda}= P^sH^n(F)_{\lambda},$$
for
$\lambda =\exp(-2 \pi ik/d)$ with 
$$0<k\leq k_0=d-\left[ \frac{d}{2} \right]-1.$$
More precisely, in these conditions one has
 $\dim Gr_F^sH^n(F)_{\lambda}= \dim Gr_P^sH^n(F)_{\lambda}= \dim M(f)_{n_1d+k-n-1}=\dim M(h)_{n_1d+k-n-1},$
 where $h$ is a homogeneous polynomial in $S$ of degree $d$ such that the hypersurface $D_h:h=0$
 is smooth.

\end{lem}

\proof In view of the inclusions in \eqref{HP1} and using Proposition \ref{propA}, it is enough to establish the inequality
$$\dim Gr_F^sH^n(F)_{\lambda} \geq \dim Gr_P^sH^n(F)_{\lambda}.$$
We know that
$$\dim Gr_P^sH^n(F)_{\lambda}=\dim E_{\infty}^{n_1+2,n_1}(f,k) \leq \dim E_1^{n_1+2,n_1}(f,k)=\dim H^{n+1}(K^*(f))_{n_1d+k}=$$
$$=\dim M(f)_{n_1d+k-n-1}= \dim M(h)_{n_1d+k-n-1}$$
where $ h \in S_d$ denotes a polynomial such that the associated hypersurface $D_h:h=0$ is smooth.
Indeed, the last equality follows from Corollary 2.2. $(i)$ in \cite{DSt3}, which is a direct consequence of Theorem \ref{thmA}.

On the other hand, using an argument similar to Proposition 4.1 in \cite{D3}, which goes back to Lemma 3.6. in \cite{Sa2}, it follows that for $\be \ne \pm 1$, the mixed Hodge structure induced on 
\begin{equation} 
\label{mhs1}
H^n(F)_{\be, \overline \be}=H^n(F)_{\be} \oplus H^n(F)_{ \overline \be}
\end{equation}
is pure of weight $n$. Therefore
$$\dim Gr_F^sH^n(F)_{\lambda}=h^{n_1+1,n_1}(H^n(F)_{\lambda})=h^{n_1+1,n_1}(H^n(\tilde D)_{\lambda})$$
in view of Corollary 1.2 in \cite{DL}. To compute the last equivariant Hodge number we use Proposition 5.2 in  \cite{DL}. The first term in the sum giving $h^{n_1+1,n_1}(H^n(\tilde D)_{\lambda})$ is the
corresponding number computed for the smooth hypersurface $\tilde D_h: \tilde h=h(x)-t^d=0$.
Using the standard identification going back to Griffiths \cite{Gr},
$$H^{n_1+1,n_1}(H^n(\tilde D_h))=H^{n+2}(K^*(\tilde h))_{(n_1+1)d},$$
and recalling that taking the $\lambda$-eigenspace means to look at forms of the form
$$g(x)t^{d-k-1}dx_0\wedge...\wedge dx_n\wedge dt,$$
where $g(x)$ is homogeneous of degree $n_1d+k-n-1$, we get
$$h^{n_1+1,n_1}(H^n(\tilde D_h)_{\lambda})= \dim M(h)_{n_1d+k-n-1}.$$ 
We show now that the other terms in the sum giving $h^{n_1+1,n_1}(H^n(\tilde D)_{\lambda})$ are trivial.
These terms are of two types:

(a) $h^{p,q}(H^{n+1}(\tilde D)_{\lambda})$, which are zero  
since a  $\be$-eigenspace  of the group $ H^{n+1}(\tilde D)$ under the $\mu_d$-action maybe nontrivial only for $\be=\pm 1$, see Example 6.3.24 in \cite{D1}, Theorem 1.1 in \cite{DL} and Theorem 4.1 in \cite{DSt3}. 

(b) $h^{p,q}(H^n(F(d))_{\lambda})$, where $F(d)$ is the affine Milnor fiber given by
$$g(y,t)=y_1^2+...y_n^2+t^d-1=0$$
in $\C^{n+1}$ with the corresponding $\mu_d$-action, i.e. 
$$\be (y_1,...,y_n,t)=(y_1,...,y_n,\be t)$$
replacing the monodromy action when eigenspaces are considered.
Since this is the Milnor fiber of an isolated weighted homogeneous singularity whose link is a rational homology sphere, we have $h^{p,q}(H^n(F(d))_{\lambda})=0$ for $p+q \ne n$. It follows that we have just one
such number to investigate, namely  $h^{n_1+1,n_1}(H^n(F(d))_{\lambda})$. Using the weights $wt(y_j)=2$ and $wt(t)=d$, we get using \cite{St}
$$h^{n_1+1,n_1}(H^n(F(d))=\sum_{j=1,2d-1}\dim M(g(y,t))_{j-d-2}.$$
The $\lambda$-eigenspace should come from the monomial $t^{d-k-1}$, of degree $2d-2k-2$.
Our condition on $k$ implies that $2d-2k-2>d-3$, hence this monomial is not giving a contribution to $h^{n_1+1,n_1}(H^n(F(d))$, i.e.
$h^{n_1+1,n_1}(H^n(F(d))_{\lambda})=0$. Moreover $2k<d$ in order to avoid the case $\lambda=-1$.

This shows that $\dim Gr_F^sH^n(F)_{\lambda}=\dim M(h)_{n_1d+k-n-1}$, completing the proof of Lemma \ref{L}.

\endproof

Now we give the proof of Theorem \ref{thmB}. Lemma \ref{L} implies that
$$E_1^{n_1+2,n_1}(f,k_0)=H^{n+1}(K^*(f))_{n_1d+k_0}=E_{\infty}^{n_1+2,n_1}(f,k_0).$$
Moreover, $E_1^{n_1+2+e,n_1-e}(f,k_0)=E_{\infty}^{n_1+1+e,n_1-e}(f,k_0)$ for all $e=1,2,...,n_1$
by similar (and simpler) computations based on Proposition \ref{propA}. It follows that all the differentials in the spectral sequence $E_r(f,k_0)$ starting from $E_r^{n_1+1,n_1}(f,k_0)$ are $0$.
Since $E_{\infty}^{n_1+1,n_1}(f,k_0)=0$ as well (the only eigenvalues of the monodromy on $H^{n-1}(F)$ are $\pm 1$), we get
$$E_1^{n_1+1,n_1}(f,k_0)=H^n(K^*(f))_{n_1d+k_0}=0$$
which completes the proof of Theorem \ref{thmB} . Indeed, recall that if the coordinates $x_0,...,x_n$ are choosen such that the hyperplane at infinity $H_0:x_0=0$ is transversal to $D$, then the multiplication by $x_0$ induces an injection
$H^n(K^*(f))_{s-1} \to H^n(K^*(f))_{s}$ for any $s$ (the dual statement for the homology is part of Corollary 11
in \cite{CD}).

Before proceeding, we recall the following notions, introduced in \cite{DSt2}.
\begin{definition}
\label{def}

For a hypersurface $D:f=0$ of degree $d$ with isolated singularities we introduce three integers, as follows:

\noindent (i) the {\it coincidence threshold} $ct(D)$ defined as
$$ct(D)=\max \{q~~:~~\dim M(f)_k=\dim M(h)_k \text{ for all } k \leq q\},$$
with $h$  a homogeneous polynomial in $S$ of degree $d=\deg f$ such that $D_h:h=0$ is a smooth hypersurface in $\PP^n$.

\noindent (ii) the {\it minimal degree of a nontrivial syzygy} $mdr(D)$ defined as
$$mdr(D)=\min \{q~~:~~ H^n(K^*(f))_{q+n}\ne 0\}$$
where $K^*(f)$ is the Koszul complex of $f_0,...,f_n$ with the natural grading defined in \cite{DSt2}.

\end{definition}

 Moreover
it is easy to see that one has 
\begin{equation} 
\label{REL}
ct(D)=mdr(D)+d-2.
\end{equation} 
In practice, for a given polynomial $f$, it is easy to compute $ct(D)$ using a number of computer algebra softwares.

\begin{cor}
\label{corA}
Let $D:f=0$ be a nodal hypersurface of degree $d>2$ in $\PP^n$.
If $n=2n_1+1$ is odd, then  
$$ct(D) \geq (n_1+2)d-\left[ \frac{d}{2} \right]-n-2.$$
\end{cor}

\begin{ex}
\label{ex1}
Let $\CC(3,d)$ be the Chebyshev surface of degree $d$ in $\PP^3$  as defined in \cite{DSt3}.
Then for $3 \leq d \leq 20$, numerical computation shows that one has
$$ct(\CC(3,d))=3d-\left[ \frac{d}{2} \right]-5,$$
i.e. we have equalities for these cases in Corollary \ref{corA}. It follows that in any such case the bound for the vanishing in Theorem \ref{thmB} is sharp, namely
$$H^n(K^*(f))_{m+1} \ne 0$$
for  $m =(n_1+1)d-\left[ \frac{d}{2} \right]-1.$
\end{ex}

\begin{ex}
\label{ex1bis}
Let $D$ be a nodal hypersurface of degree $d$ in $\PP^n$  having exactly one singularity.
Then it is known that
$$mdr(f)=n(d-2) \text{ and } ct(D)=(n+1)(d-2),$$
see \cite{DSt2}, i.e. a substantially bigger number that the lower bound given in Corollary \ref{corA}. A better lower bound for $ct(D)$ than that given by Corollary \ref{corA} in the case of a nodal hypersurface with not too many nodes is given in \cite{DAG}. More precisely, it is shown that 
$$mdr(f) \geq n(d-2)+1 - \sharp A_1  \text{ and } ct(D) \geq (n+1)(d-2)+1 - \sharp A_1,$$
where $\sharp A_1$ denotes the number of nodes of the nodal hypersurface $D$. Note that these inequalities are sharp when $ \sharp A_1=1$. For the case of several nodes which are linearly independent, see \cite[Proposition 1]{DS4}.
\end{ex}

\begin{rk} \label{sernesi}

Let $D:f=0$ be a degree $d$ hypersurface in $\PP^n$ having only isolated singularities. Let $\hat J$ be the saturation of the Jacobian ideal $J$ of $f$. Then 
the vector space $\hat J_d/J_d$ is naturally identified with the space of first order locally trivial deformations of $D$ in $\PP^n$ modulo those arising from the above $PGL(n+1)$-action, see E. Sernesi \cite{Se}.
 The dimension of the vector space $\hat J_d/J_d$ can be determined as follows.
$$\dim \hat J_d/J_d=\dim M(f)_{d}-\dim M(h)_d+\dim M(f)_{T-d}-\tau(D),$$
where $h$ is as in Definition \ref{def}, see G. Sticlaru \cite{Sti} for this formula and a number of interesting examples of rigid and non-rigid  hypersurfaces. Here $\tau(D)$ is the total Tjurina number of the hypersurface $D$, e.g. the number of nodes for a nodal hypersurface.

\end{rk} 

\section{Hodge theory of nodal hypersurfaces and the proof of Corollary \ref{corSeveri} and Theorem \ref{corNobile}} \label{sec5}

Let $\I \subset \OO_{\PP^n}$ be the reduced ideal sheaf of the set of nodes $\NN=Sing D \subset \PP^n$. Set $I_k{(i)}=H^0(\PP^n,\I^i(k))\subset S_k$ and define $I^{(i)}=\oplus_kI_k{(i)}$,
a homogeneous ideal in the polynomial ring $S$.

For a degree $d$ nodal hypersurface $D$ in $\PP^n$, one of the main results in \cite{DSW} describe the graded pieces $Gr_F^pH^n(U)$ of the top cohomology of the complement $U=\PP^n \setminus D$ with respect to the Hodge filtration $F$ in terms of purely algebraic objects, namely one has
\begin{equation} 
\label{GR}
Gr_F^pH^n(U,\C)=(I^{(q-m+1)}/I^{(q-m)}J_f)_{(q+1)d-n-1}
\end{equation} 
for $q=n-p>m:=\left[ \frac{n}{2} \right]$ under a certain condition (B), see Theorem 2 in \cite{DSW}.

Recall that for a finite set of points $\NN \subset \PP^n$
we denote by 
$$\defect S_m(\NN)=|\NN| - \codim \{h \in S_m~~|~~ h(a)=0 \text{ for any } a \in \NN\},$$ 
the {\it defect (or superabundance) of the linear system of polynomials in $S_m$ vanishing at the points in $\NN$}, see \cite{D1}, p. 207. This positive integer is called the {\it failure of $\NN$ to impose independent conditions on homogeneous polynomials of degree $m$} in \cite{EGH}. 
There is a close relationship between defects $\defect S_m(\NN)$ and the syzygies described by $H^n(K^*(f))$, see \cite{DSt2} for nodal hypersurfaces and \cite{D4} for projective hypersurfaces with isolated singularities. More precisely, for the nodal hypersurfaces one has the following.
\begin{equation} 
\label{defect}
 \defect S_k(\NN)= \dim H^n(K^*(f))_{nd-n-1-k},
\end{equation} 
for $0 \leq k \leq nd-2n-1$ and $ \defect S_k(\NN)=0$ for $k > nd-2n-1$.

The discussion following the statement of Theorem 2 in \cite{DSW} shows that in fact (B) is equivalent to the condition
\begin{equation} 
\label{B}
(B'): \defect S_e(\NN)=0,
\end{equation} 
where $e=\left[ \frac{n}{2} \right](d-1)-p$ and $\NN $ is the set of nodes of $D$.

One has the following consequence of Theorems \ref{thmA} and \ref{thmB}.

\begin{thm}
\label{thmC}
Let $D:f=0$ be a nodal hypersurface of degree $d>2$ in $\PP^n$ and assume that
$q=n-p>m:=\left[ \frac{n}{2} \right]$.
\smallskip

\noindent (i) If $n=2n_1$ is even, then the isomorphism \eqref{GR} always hold.

\smallskip

\noindent (ii) If $n=2n_1+1$ is even, then the isomorphism \eqref{GR} holds if either
$$p\leq n -\left[ \frac{n}{2} \right] -\left[ \frac{d}{2} \right] ,  $$
or 
$$\sharp A_1 \leq (n_1+2)(d-1)+p-n.$$
In particular, the isomorphism \eqref{GR} holds always for $d=3$, and for $d=4$ with the additional condition $\sharp A_1 \leq n_1+p+5$.

\end{thm}

\proof When $n=2n_1$ is even, then $p<n_1$ and hence $e=n_1(d-1)-p=n_1d-n_1-p>n_1d-n.$
Use now Corollary 2.2 (ii) in \cite{DSt3} which says that
$\defect S_k(\NN)=0$ for $k\geq n_1d-n$.

\smallskip

When $n=2n_1+1$, we have $e=n_1(d-1)-p$.
On the other hand, one know that $\defect S_k(\NN)=0$ if and only if 
$k\geq T-ct(D)$, see Theorem 1.5 in  \cite{DSt2}.  Corollary \ref{corA} implies
$$T-ct(D) \leq 2(n_1+1)(d-2)-\left( (n_1+2)d-\left[ \frac{d}{2} \right]-n-2\right)=n_1d-n+\left[ \frac{d}{2} \right].$$
Hence $e \geq n_1d-n+\left[ \frac{d}{2} \right]$ as soon as $p\leq n -\left[ \frac{n}{2} \right] -\left[ \frac{d}{2} \right] .  $ For the claims involving $\sharp A_1$, use Example \ref{ex1bis}.

\endproof

Note that Example 4.7 in \cite{DSW} shows that (B) may not be satisfied by a nodal surface, where $n=3$,  $d=4$, $p=1$ and  $\sharp A_1$ is large. This shows that the first condition in Theorem \ref{thmC} (ii) is sharp,
as in the case at hand we get $1 \leq 3-1-2$, which fails just by $1$.
The following result generalizes Corollary 1 in \cite{DSW}.
\begin{cor}
\label{corAbis}
Let $D:f=0$ be a nodal hypersurface of degree $d \geq 5$ in $\PP^n$.
If $n$ is odd and $\sharp A_1 \leq n+p+6$, then the isomorphism \eqref{GR} holds.
\end{cor}

Now we give the proof of Corollary \ref{corSeveri}. When $n=3$, then by Theorem \ref{thmB} we get $H^3(K^*(f))_m=0$ for any $m \leq 2d -\left[ \frac{d}{2} \right]-1$. Using formula \eqref{defect}, this implies $\defect S_k(\NN)=0$ for $k \geq d +\left[ \frac{d}{2} \right]-3$, which is exactly the claim in Corollary \ref{corSeveri}.
Note that for $d >2$ we have 
$$d +\left[ \frac{d}{2} \right]-3 \leq 2d-5,$$
with equality only for $d=3$ and $d=4$. In other words, our result is stronger than Severi's for any $d\geq 5$.

The proof of the first two claims in Theorem \ref{corNobile} follows from Corollary \ref{corSeveri} and the discussion in \cite[(3.6)]{Nobile}. Indeed, one can check that the inequality
$$d +\left[ \frac{d}{2} \right]-3 \leq d$$
holds exactly for $d \leq 7$. To prove the third claim, note that one has
$$\dim m_{R(X)}/m_{R(X)}^2-\dim R(X)=\defect S_d(\NN),$$
for $d \geq 5$, see \cite[Theorem 4.2]{BW}, where  $m_{R(X)}$ denotes the maximal ideal of the  complete local ring $R(X)$.

\begin{rk} \label{Krulldim}

\begin{enumerate} 

\item Note  that the Krull dimension $\dim R(X)$ of the formal moduli space $R(X)$ of $X$
is given by 
$N-15$, see \cite[Corollary 2.11]{BW}, which is in general different from $\dim V_{d,n}$.
In other words, the (embedded) deformation theory of the nodal surface $D$, reflected in the germ of $V_{d,n}$ at $D$, is quite different from the deformation theory of the surface $X$, reflected in the local ring $R(X)$.

\item For a surface $D$ in $\PP^3$ of degree $d \leq 4$ with only rational double points, one has that the corresponding minimal surface $X$ is unobstructed, see \cite[Example (4.7)]{BW}. In the same Example, the authors produce a quintic $D$ with 10 $A_4$-singularities
which is obstructed. 

\item It is very likely that the  second claim in Theorem \ref{corNobile} can be upgraded to the inclusion $V_{d,n} \subset \overline V_{d,n'}$, as it is the case for the Severi varieties of nodal plane curves, see \cite{Nobile0} or \cite[Theorem 4.7.18]{Se0}. To do this, one should check that the methods described in \cite{Nobile0} or \cite{Se0.5} work in the case of nodal surfaces as well.

\end{enumerate} 
\end{rk}

\end{document}